# A GENERALIZATION OF THE
# "BROUWER – SCHAUDER – TYCHONOFF" FIXED-POINT THEOREM

### RANJIT VOHRA

**ABSTRACT**. We prove a new fixed – point result for the image $\mathrm{Im}(\varphi)$ of any continuous function $\varphi$ from K to $K \times K$, where K is a compact convex subset of a Hausdorff locally convex space, provided that the projection of $\mathrm{Im}(\varphi)$ to the first factor is onto, and a condition on the convex hull of $\mathrm{Im}(\varphi)$ holds. A special case of our result is the Brouwer-Schauder-Tychonoff fixed-point theorem for continuous functions $f : K \to K$.

## §**1**. INTRODUCTION

We know, from the Brouwer-Schauder-Tychonoff theorem (see below for theorem statement), that any continuous function $f : K \to K$, where K is a compact convex subset of a Hausdorff locally convex space, admits a fixed point, meaning that the graph of $f$ intersects the diagonal $\Delta \subset K \times K$. Intuitively, it seems plausible that any squiggly continuous image of K in $K \times K$, regardless of whether or not this image describes a function $f : K \to K$ in the target $K \times K$, should also intersect $\Delta$. We make this intuition precise in our Theorem I below. Furthermore, since the Brouwer-Schauder-Tychonoff theorem can be regarded as a corollary of the Kakutani-Fan-Glicksberg ("KFG") fixed-point theorem for correspondences, we note (see Remark #2 in §3 below) that the hypotheses of our Theorem I are not comparable to those of KFG.

Subject Classification: 46A99, 54A99, 91B50

## §**2**. REMINDERS

We recall the following two results.

BROUWER – SCHAUDER – TYCHONOFF FIXED-POINT THEOREM. Suppose:

- X is a Hausdorff locally convex space;
- K is a nonempty compact, convex subset of X;
- $f : K \to K$ is a continuous function.

Then $f$ has a fixed point.

*Proof.* See original papers [4], [11], [12]; and also 17.56, p.583 in [1].

KAKUTANI – FAN – GLICKSBERG FIXED-POINT THEOREM. Suppose:

- X is a Hausdorff locally convex space;
- K is a nonempty compact, convex subset of X;
- $\Gamma : K \twoheadrightarrow K$ is a correspondence that has closed graph and nonempty, convex values.

Then $\Gamma$ has a fixed point.

*Proof.* See original papers [10], [6], [7]; and also 17.55, p.583 in [1]

FACT: Every locally convex space (even if not Hausdorff) is locally connected, since each point has a local base of convex (hence connected) neighborhoods. In any locally connected space, each component of any open set is open. See e.g. p.199*ff* in [13]

## §3.                      THE MAIN RESULT

***Remark*** #1. In the sequel we will see the following situation. There is a given function $\varphi : K \to K \times K, \ x \mapsto (\varphi_1(x), \varphi_2(x))$ such that $\varphi_1$ is onto. The image of $\varphi$, denoted $\text{Im}(\varphi)$, is a subset of $K \times K$ and can be regarded as the graph of the correspondence $\Psi : K \twoheadrightarrow K,$ given by $t \mapsto \varphi_2(\{\varphi_1^{-1}(t)\})$. Thus $\text{Im}(\varphi) = \text{gr}(\Psi)$; where the correspondence $\Psi$ is nonempty – valued (using that $\varphi_1$ is onto).

THEOREM **I**. Suppose:

- X is a Hausdorff locally convex space;
- K is a nonempty compact, convex subset of X;
- $\varphi : K \to K \times K, \ x \mapsto (\varphi_1(x), \varphi_2(x))$ is a continuous function such that $\varphi_1$ is onto; & the correspondence $\Psi : K \twoheadrightarrow K, \ t \mapsto \varphi_2(\{\varphi_1^{-1}(t)\})$ has the following property: for each $t$ in the domain of $\Psi$, $t \notin \Psi(t)$ implies $t \notin \text{co}(\Psi(t))$.

Then $\Psi$ has a fixed point.

*Proof.* By Remark #1, the correspondence $\Psi : K \twoheadrightarrow K, \ t \mapsto \varphi_2(\{\varphi_1^{-1}(t)\})$ is nonempty - valued. Note also that $\text{gr}(\Psi)$ is connected, since $\text{gr}(\Psi) = \text{Im}(\varphi)$, and $\text{Im}(\varphi)$ is connected because it is the continuous image of the convex (hence connected) set K.

Now suppose, toward a contradiction, that the theorem hypotheses hold, and that $\Psi$ has no fixed point, meaning that $\mathrm{gr}(\Psi) \cap \Delta = \varnothing$, where $\Delta$ denotes the diagonal in the compact product $K \times K$. Then $t \notin \Psi(t)$ and indeed $t \notin \mathrm{co}(\Psi(t))$ holds for all $t$ in the domain K of $\Psi$, using the hypotheses on $\Psi$.

Put $\quad \Psi'(t) := \mathrm{co}(\Psi(t))$, and define a new correspondence $\Psi': K \twoheadrightarrow K$ by $\Psi' := \bigcup_{t \in K} \mathrm{co}(\Psi(t))$. We see that $\mathrm{gr}(\Psi')$ has two important properties:

(P1)  $\mathrm{gr}(\Psi')$ is connected, since it is the union of the connected set $\mathrm{gr}(\Psi)$ and the convex sets $\mathrm{co}(\Psi(t))$ that are attached to points of $\mathrm{gr}(\Psi)$.

(P2)  $\mathrm{gr}(\Psi') \cap \Delta = \varnothing$ since $t \notin \Psi'(t)$ holds for each $t \in K$, as noted above.

Since $\Delta$ is closed in $K \times K$ (using the Hausdorff hypothesis), we see from (P2) that $\mathrm{gr}(\Psi')$ must belong to the open set $(K \times K) \backslash \Delta$. Moreover, since $\mathrm{gr}(\Psi')$ is connected (P1), we know (see Fact in §2) there is an open component $\mathcal{U}$ of $(K \times K) \backslash \Delta$ such that $\mathrm{gr}(\Psi') \subseteq \mathcal{U} \subseteq (K \times K) \backslash \Delta$.

The idea is to "fatten" $\Psi'$ by adding open sets to its boundary $\partial\Psi'$, in order to produce a new correspondence $\Psi^*$ that contains $\Psi'$ and has open graph. Precisely, to each $z$ of $\partial\Psi'$, we attach a basic open disk $B_z$ (of the locally convex space $X \times X$). Thus for each $t$ in K, and using notation $(B_z \mid t)$ for the restriction $B_z \cap (\{t\} \times K)$, we obtain:

$$\Psi^*(t) := \Psi'(t) \cup \left( \bigcup_{z \in \partial(\Psi'(t))} (B_z \mid t) \right) \qquad (\dagger)$$

where each $B_z$ is a basic open, convex disk of $X \times X$ that is centered at point $z$ on the boundary of $\Psi'$, and is small enough that $\Psi^*(t) \subseteq \mathcal{U} \subseteq (K \times K) \backslash \Delta$.

We claim that $\underline{t \notin \mathrm{co}(\Psi^*(t))}$ for any $t \in K$. To see this, note that

$$\mathrm{co}(\Psi^*(t)) = \mathrm{co}\left( \Psi'(t) \cup \left( \bigcup_{z \in \partial(\Psi'(t))} (B_z \mid t) \right) \right) \qquad \text{using } (\dagger)$$

$$= \Psi'(t) \cup \left( \bigcup_{z \in \partial(\Psi'(t))} (B_z \mid t) \right)$$

as the expression on the last line is a convex set that (by choice of the disks $B_z$) does not contain $t$.

Put $\Psi^* := \bigcup_{t \in K} \Psi^*(t)$. Evidently $\Psi^*$ is a nonempty – valued correspondence $K \twoheadrightarrow K$ with graph $\mathrm{gr}(\Psi^*)$ that is open and connected, by construction. Thus we obtain the inclusions $\mathrm{gr}(\Psi^*) \subseteq \mathcal{U} \subseteq (K \times K) \backslash \Delta$, and (as above) that $t \notin \mathrm{co}(\Psi^*(t))$ holds for each $t$ in the domain K of $\Psi$.

Let $\mathcal{H} : \mathrm{K} \twoheadrightarrow \mathrm{K}$, $y \mapsto \mathrm{K} \backslash (\Psi^*)^{-1}(y)$ be the "inverse complement correspondence"[1] of $\Psi^* : \mathrm{K} \twoheadrightarrow \mathrm{K}$. Observe that Lemma 17.47 in [1]) applies to $\Psi^*$, and since $t \notin \mathrm{co}(\Psi^*(t))$ for all $t \in \mathrm{K}$, we see $\mathcal{H}$ is a KKM correspondence. Furthermore (using that $\Psi^*$ has open graph) Theorem 17.46 of [1] implies that there is some point $x$ in the domain K of $\Psi^*$ (and hence also of $\Psi$) that is empty–valued, contradicting that $\Psi$ is nonempty–valued. We therefore conclude from this contradiction that our initial supposition $\mathrm{gr}(\Psi) \cap \Delta = \varnothing$ is untenable, and that $\Psi$ must have a fixed point. $\square$

COROLLARY. Suppose the function $\varphi_1$ is the identity map. Then $\Psi = \varphi_2$, whereby $\Psi$ is a continuous function from K to K, and we recover the Brouwer-Schauder-Tychonoff fixed point theorem (see § 2) for $\Psi$.

**Remark** #2. Although the Brouwer – Schauder – Tychonoff fixed point theorem is a special cass of both the Kakutani – Fan – Glicksberg theorem and of our Theorem I, the hypotheses of our theorem are different from those of KFG: neither implies or is implied by the other. The hypotheses of the FKG theorem require that the correspondence $\Gamma$ has convex values, whereas our Theorem I has the weaker requirement that for each $t$ in the domain of the correspondence $\Psi$, $t \notin \Psi(t)$ implies $t \notin \mathrm{co}(\Psi(t))$. On the other hand, the correspondence $\Gamma$ of the FKG theorem is not necessarily the image of a continuous function $\varphi : \mathrm{K} \to \mathrm{K} \times \mathrm{K}$, so FKG is more general in this regard.

# REFERENCES

[1] C.D. Aliprantis, K.C. Border, Infinite Dimensional Analysis, 3$^{\mathrm{rd}}$ ed., Springer, 2006.

[2] M. Balaj, S. Muresan, Generalizations of the Fan-Browder Fixed Point Theorem & Minimax Inequalities, Archivum Mathematicum, Tomus 41 (2005), 399-407.

[3] K.C. Border, Fixed point theorems with applications to economics and game theory, Cambridge University Press, 1985

[4] L.E.J. Brouwer, 1912, Uber Abbildung von Mannigfaltieiten, Mathematische Annalen 71:165 -168.

[5] K. Fan, Fixed point and minimax theorems in locally convex topological spaces, Proc. Of the National Academy of Sciences, USA 38:121-126, (1952)

[6] K. Fan, A generalization of Tychonoff's fixed point theorem, Math. Ann. 142 (1962) 305–310.

[7] I.L. Glicksberg, A further generalization of the Kakutani fixed point theorem with application to Nash equilibrium, Proc. Amer. Math. Soc. 3 (1952) 170–174.

---

[1] See e.g. p.578 in [1].

[8]   B.R. Halpern and G.M. Bergman, A Fixed Point Theorem for Inward and Outward Maps, Trans. Amer. Math. Soc., 130, 353-358 (1968).

[9]   C.J. Himmelberg, Fixed points of compact multifunctions, J. Math. Anal. Appl. 38(1972), 205–207

[10] S. Kakutani, A generalization of Brouwer's fixed point theorem, Duke Mathematical Journal, 8:457-459, 1941.

[11] J. Schauder, 1930, Der Finxpunksatz in Functionalraumen, Studia Mathematica 2; 171-180.

[12] A. Tychonoff, 1935, Ein Fixpunktsatz, Mathematische Annalen 111:767-776.

[13] S. Willard, General Topology, Addison Wesley, 1970 .

R. VOHRA

P.O. Box 240, Storrs CT 06268

rvohra@bridgew.edu